\documentclass[a4paper,twoside,12pt]{article}
\usepackage[english]{babel}
\usepackage{graphicx}
\usepackage{bbding}
\usepackage[latin1]{inputenc}
\usepackage{kpfonts}
\usepackage{fancyhdr}
\usepackage{amsfonts,amsmath,amssymb,color,epsfig}
\usepackage{mathrsfs}
\usepackage{eufrak}
\usepackage{dsfont}
\usepackage{upgreek}
\usepackage{fancyhdr}
\usepackage[dvipsnames]{xcolor}
\newtheorem{theorem}{Theorem}[section]
\newtheorem{proposition}[theorem]{Proposition}

\newtheorem{remark}[theorem]{Remark}
\newtheorem{lemma}[theorem]{Lemma}

\def\N{{\mathbb N}}

\def\R{{\mathbb R}}
\def\C{{\mathbb C}}

\def\K{{\mathbb K}}

\def\1{\mathds{1}}

\title{Product of Independent Cauchy-Lorentz Random Matrices }
\date{}
\definecolor{mon}{rgb}{10,126,140}
\begin{document}
\author{Mohamed BOUALI}
\maketitle
\begin{abstract}
We investigate the product of $n$ complex non-Hermitian, independent random matrices, each of
size $N_i\times N_{i+1}$ $(i=1,...,n)$, with independent identically distributed Cauchy entries (Cauchy-Lorentz matrices). The joint
probability distribution of the complex eigenvalues of the product matrix is found to be given by a
determinantal point process as in the case of a single Cauchy-Lorentz matrix, but with
weight given by a Meijer G-function depending on $n$ and $N_i$.
\end{abstract}
\section{Introduction}
The topic of random matrix theory (RMT) enjoys an increasing number of applications in physics,
mathematics and other sciences, and we refer to \cite{ak} and \cite{fo} for a recent compilation. This statement holds
both for problems with real eigenvalues as well as with complex eigenvalues. In the latter case the
classical ensemble was introduced by Ginibre \cite{gi} who considered complex non-Hermitian matrices $X$
of size $N\times N$, with all matrix elements having independent normal distributions. However, for certain
applications it is not sufficient to introduce a single random matrix, e.g. when considering transfer
matrices. The problem of studying such products of random matrices is as old as RMT itself and was
introduced by Furstenberg and Kesten \cite{fu}.

One problem one has to face is that the product matrix often has less symmetry than the individual
matrices. For example the product $P$ of two Hermitian matrices is in general no longer Hermitian,
and thus acquires a complex spectrum. Therefore in the literature two types of products have been
considered, those which are Hermitised by considering $P^\ast P$ with real spectra, and those which are
studying the complex eigenvalues of $P$ itself.

In this work we will consider the case of complex eigenvalues, by multiplying a fixed number $n$ of independent Cauchy-Lorentz matrices: $Y_n = X_1X_2 . . .X_n$. For such ensemble we give the limit of the one point correlation function (state density) $\rho_{n,N}(z)$ as $N$ goes to $+\infty$. One will see that such density do not require rescaling  like as in the case of product of $n$-independent Ginibre matrices where a rescaling of order $N^{\frac n2}$ is necessary. In fact for $n=1$ the density is $\displaystyle\omega_{1,N}=\frac1{(1+|z|^2)^{N+1}}$, and It can be seen within the Coulomb gas picture, that the Vandermonde
determinant is a logarithmic two-body repulsive term and the weight function is essentially a logarithmic (one-body) confining potential, that needs the explicit dependence with $N$ to be able to compensate the repulsive term. Maybe the same hold for $n>1$ and, it is rather natural that we get asymptotically the same expression for the density of
eigenvalues for all, since we have a modified confining potential for each $N$.

In section two, we give the density of eigenvalues of the product of $n$ independent Cauchy-Lorentz matrices and we express the weight $\omega_{n,N}$ as a multiple integral and then as a Meijer $G$-function which depends on the parameters $n$ and $N_i$. Moreover we investigate as $N\to+\infty$ $(N=N_1)$ the asymptotic of the weight $\omega_{n,N}(z)$ for two regimes, the first one for $z$ in every compact of $\C$ and the second one at the origin which will be need to rascal $z$ by the factor $N^{-\frac n2}$.

In section three, Using the method of orthogonal polynomials
we compute all eigenvalues density correlation functions exactly for finite $N_i$ and fixed $n$, which will be expressed as a generalized hypergeometric function ${}_nF_{n-1}$. In the large-$N$
limit at fixed $n$ and fixed $\nu_i=N_i-N$, $(N=N_1)$ we first determine the macroscopic limit of the one point correlation function (global density of eigenvalues) and the microscopic correlation functions in the bulk, and at the origin. The case of bulk limit are identical to that of the Ginibre ensemble with $n = 1$ and
thus universal. However the microscopic correlations we find at the origin differ for each $n > 1$
and generalise the known Bessel Kernel in the complex plane for $n = 2$ to a new hypergeometric kernel
${}_0F_{n-1}$.
\section{Density of eigenvalues}
Let $N\in\N$, and $\varphi_N:M_{N\times m}(\Bbb C)\longrightarrow\R$ be a function. One consider the following density
$$P_N(X)=\varphi_N(XX^*)dXdX^*,$$
on the space of $N\times m$, complex rectangular matrices. Assume that $\varphi_N$ is unitary invariant and for $n\in\N$, consider $X_1,...,X_n$ as complex rectangular matrices of orders $N_i\times N_{i+1}$, $(i=1,...,n)$, then it has been proved in \cite{ak1} and \cite{ip} and others references that, after the  assumption $N=N_1\leq...\leq N_{n+1}$, the density of the complex eigenvalues of the matrices $Y_n=X_1\cdots X_n$ is given by
$${\cal P}_{n,N}(z_1,...,z_N)=\prod_{j=1}^Nd^2z_j\omega_{n,N}(z_j)\prod_{1\leq i<j\leq N}|z_i-z_j|^2,$$
where $z_1,...,z_N$ are the complex eigenvalues of the matrix $Y_n$ and
$$\omega_{n,N}(z)=\int_{\C^n}\prod_{k=1}^n\varphi_N(|x_k|^2)\delta(z-\prod_{k=1}^nx_k)d^2x_1...d^2x_n.$$

For $i=1,...,n$, Let $X_i$ be a matrix in the class of $N_i\times N_{i+1}$ complex random Cauchy-Lorentez matrices with density
 $$p_{N_i}(X_i)=\frac1{C_{n,N_i}}\frac1{\det(1+XX^*)^{N_i+1}}dX_idX^*_i,$$
where $\displaystyle C_{n,N_i}= \frac{2^{(N_i+1)N_in}}{(\pi\Gamma(N_i+1))^{n}}.$ We assume that $N_1\leq N_2\leq...\leq N_{n+1}$ and we
denote in the sequel $N=N_1$, $\nu_i=N_i-N$ and $\displaystyle\alpha_n=\sum_{i=1}^n\nu_i$ (remark that $\nu_1=0$).

The pdf density of complex eigenvalues of the product random independent matrices $Y_n=X_1\cdots X_n$, is
$${\cal P}_{n,N}(z_1,\cdots z_N)=\frac1{(C_{n,N})^n}\prod_{i=1}^Nd^2z_i\omega_{n,N}(z_i)\prod_{1\leq i<j\leq N}|z_i-z_j|^2,$$
with $z_1,...,z_N$ the complex eigenvalues of the matrix $Y_n$,
and  $$ \omega_{n,N}(z)=\int_{\mathbb C^n}d^2x_1...d^2x_n\prod_{j=1}^n\frac{|x_i|^{2\nu_j}}{\big(1+|x_j|^2\big)^{\nu_i+N+1}}\delta\Big(z-\prod_{j=1}^nx_j\Big).$$
The moments of the density $\omega_n(z)$ are given for every $k\leq N-1$ by
$$m_{2k}=\int_\C|z|^{2k}\omega_n(z)d^2z=\prod_{i=1}^n\int_{\C}\frac{|x_i|^{2k+2\nu_i}}{\big(1+|x_i|^2\big)^{\nu_i+N+1}}d^2x_i,$$
$$m_{2k}=\pi^n\big(\Gamma(N-k)\big)^n\prod_{i=1}^n\frac{\Gamma(\nu_i+k+1)}{\Gamma(\nu_i+N+1)}.$$
For $k\in\N$, the $k-$point correlation function is defined by
$$\rho_{n,k,N}(z_1,\cdots ,z_k)=\frac{N!}{(N-k)!}\frac1{C_{n,N}}\int_{\C^{N-k}}{\cal P}_{n,N}(z_1,\cdots ,z_N)d^2z_{k+1}\cdots d^2z_N=\det_{1\leq i,j\leq k}\Big(K_{n,N}(z_i,z_j)\Big).$$
The correlation function is a determinantal point process on $\C$ with kernel $K_{n,N}(z_i,z_j)$, that is
$$\rho_{n,k,N}(z_1,\cdots ,z_k)=\det_{1\leq i,j\leq k}\Big(K_{n,N}(z_i,z_j)\Big),$$
where $$K_{n,N}(z_i,z_j)=\sqrt{\omega_{n,N}(z_i)\omega_{n,N}(z_j^*)}\sum_{k=1}^{N-1}p_k(z_i)p_k(z_j^*),$$
and $p_k(z)$ is the family of orthogonal polynomials with respect to the weight $\omega_{n,N}$.

Using the Mellin transform one can derive the following formula for the weight $\omega_{n,N}(z)$
\begin{proposition} For every $z\in\C$
$$\omega_{n,N}(z)=\frac{\pi^{n-1}}{\prod_{i=1}^n\Gamma(\nu_i+N+1)}G_{n,n}^{n,n}\Big(\begin{aligned}&-N,...,-N\\&0,\nu_2,\nu_3,....,\nu_n\end{aligned}\Big|\,|z|^2\Big).$$
where $G_{p,q}^{n,m}\Big(\cdots\Big|\,x\Big)$ denotes the Meijer's $G$-function.
\end{proposition}
For $n=1$, one gets $$\omega_{1,N}(z)=\frac1{\Gamma(1+N)}G_{1,1}^{1,1}\Big(\begin{aligned}&-N\\&\;\;\; 0\end{aligned}\Big|\,|z|^2\Big)=\frac1{(1+|z|^2)^{N+1}}.$$
which coincide with the density of Cauchy-Lorentz random matrix ensemble.
\section{Integral expression of the weight}
Let $r=|z|$. If we use the change of coordinates $z=re^{i\theta}$, $\displaystyle x_n=\frac{z}{x_1\cdots x_{n-1}}$ and $r_j=|x_j|$, one gets
$$\int_\C\omega_{n,N}(z)d^2z=\int_\C\prod_{j=1}^n\frac{|x_j|^{2\nu_j}}{(1+|x_j|^2)^{\nu_j+N+1}}d^2x_j.$$
Hence $$\int_\C\omega_{n,N}(z)d^2z=(2\pi)^{n-1}\int_0^{2\pi}d\theta\int_0^{+\infty}\frac{rdrr_1dr_1\cdots r_{n-1}dr_{n-1}}{r_1^2\cdots r_{n-1}^2}\prod_{j=1}^{n-1}\frac{r_j^{2\nu_j-2\nu_n}}{(1+r^2_j)^{\nu_j+N+1}}\frac{r^{2\nu_n}}{(1+\frac{r^2}{r_1^2\cdots r_{n-1}^2})^{\nu_n+N+1}}.$$
Which gives the following expression for the weight $\omega_{n,N}$,
\begin{equation}\label{eq4}\omega_{n,N}(z)=(2\pi)^{n-1}\int_0^{+\infty}dr_1\cdots dr_{n-1}\frac{|z|^{2\nu_n}}{(1+\frac{|z|^2}{r_1^2\cdots r_{n-1}^2})^{\nu_n+N+1}}\prod_{j=1}^{n-1}\frac{r_j^{2\nu_j-2\nu_n-1}}{(1+r^2_j)^{\nu_j+N+1}}.\end{equation}

Let $\varphi_{n,N}(x)=\omega_{n,N}(\sqrt x)$, with $x=|z|$, and $M_n$ the Mellin transform of $\varphi_{n,N}$. Since
$$M_n(s)=\int_0^{+\infty}\varphi_n(x)x^{s-1}dx,$$
and $$M_n(s)=(2\pi)^{n-1}\int_0^{+\infty}\int_0^{+\infty}\frac{x^{2\nu_n}}{(1+\frac{x^2}{r_1^2\cdots r_{n-1}^2})^{\nu_n+N+1}}\prod_{j=1}^{n-1}\frac{r_j^{2\nu_j-2\nu_n-1}}{(1+r^2_j)^{\nu_j+N+1}}x^{s-1}dr_1\cdots dr_{n-1}dx,$$
Thus, by the change of variable $\displaystyle u=\frac x{r_1...r_{n-1}}$, it follows that
$$M_n(s)=(2\pi)^{n-1}\int_0^{+\infty}\int_0^{+\infty}\frac{u^{2\nu_n+s-1}}{(1+u^2)^{\nu_n+N+1}}\prod_{j=1}^{n-1}
\frac{r_j^{2\nu_j+s-1}}{(1+r^2_j)^{\nu_j+N+1}}dr_1\cdots dr_{n-1}du.$$
hence
$$M_n(s)=(2\pi)^{n-1}\prod_{j=1}^{n}\int_0^{+\infty}
\frac{r_j^{2\nu_j+s-1}}{(1+r^2_j)^{\nu_j+N+1}}dr_j.$$
And $$M_n(s)=\pi^{n-1}\prod_{j=1}^{n}\frac{\Gamma(N+1-\frac s2)\Gamma(\nu_j+\frac s2)}{\Gamma(\nu_j+N+1)},$$
Moreover, by the inverse Mellin transform one has
$$\varphi_{n,N}(x)=\pi^{n-1}\frac1{2i\pi}\int_{c-i\infty}^{c+i\infty}(\Gamma(N-\frac s2))^n\prod_{j=1}^{n}\frac{\Gamma(\nu_j+\frac s2)}{\Gamma(\nu_j+N+1)}x^{-s}ds.$$
The last integral is nothing than the Meijer G-function $$\frac{\pi^{n-1}}{\prod_{i=1}^n\Gamma(\nu_i+N+1)}G_{n,n}^{n,n}\Big(\begin{aligned}
&-N,...,-N\\&0,\nu_2,\nu_3,....,\nu_n\end{aligned}\Big|\,x^2\Big).$$
\begin{lemma}
For every $z$ in compact subset of $\Bbb C$, and for fixed $\nu_2,...,\nu_n$, one has as $N\to +\infty$
$$\omega_{n,N}\big(\frac z{N^{\frac n2}}\big)\sim \frac{(2\pi)^{n-1}}{N^{\alpha_n}}|z|^{2\nu_n}\int_0^{+\infty}e^{-\sum_{j=1}^{n-1}r_j^2-\frac{|z|^2}{r_1^2\cdots r_{n-1}^2}}\prod_{j=1}^{n-1}r_j^{2\nu_j-2\nu_n-1}dr_1\cdots dr_{n-1},$$
where $\displaystyle \alpha_n=\sum_{i=1}^n\nu_i$
\end{lemma}
{\bf Proof.} From equation (\ref{eq4}), we have
$$\omega_{n,N}(\frac z{N^\frac n2})=\frac{(2\pi)^{n-1}}{N^{n\nu_n}}\int_0^{+\infty}dr_1\cdots dr_{n-1}\frac{|z|^{2\nu_n}}{(1+\frac{|z|^2}{N^{n}r_1^2\cdots r_{n-1}^2})^{\nu_n+N+1}}\prod_{j=1}^{n-1}\frac{r_j^{2\nu_j-2\nu_n-1}}{(1+r^2_j)^{\nu_j+N+1}}.$$
By the change of variables $u_i=\sqrt Nr_i$, it follows that
$$\omega_{n,N}(\frac z{N^\frac n2})=\frac{(2\pi)^{n-1}}{N^{\alpha_n}}\int_0^{+\infty}du_1\cdots du_{n-1}\frac{|z|^{2\nu_n}}{(1+\frac{|z|^2}{N\;u_1^2\cdots u_{n-1}^2})^{\nu_n+N+1}}\prod_{j=1}^{n-1}\frac{u_j^{2\nu_j-2\nu_n-1}}{(1+\frac{u^2_j}{N})^{\nu_j+N+1}}.$$
Since for fixed $\nu_j$ and for $z\in\C$, $\lim_{N\to\infty}(1+\frac{u^2_j}{N})^{\nu_j+N+1}=e^{u_j^2}$ and
 $\lim_{N\to\infty}(1+\frac{|z|^2}{N\;u_1^2\cdots u_{n-1}^2})^{\nu_n+N+1}=e^{\frac{|z|^2}{u_1^2\cdots u_{n-1}^2}}.$ Which gives the desired result by a simple application of the dominate convergence theorem on every compact set of $\C$.

In the following proposition we derive an asymptotic expression for the weight $\omega_{n,N}$ near the origin.
\begin{proposition}\label{pr01} For every $z$ in compact subset of $\C$ and for $N$ large enough
$$\omega_{n,N}\big(\frac z{N^{\frac n2}}\big)\sim \frac{\pi^{n-1}}{N^{\alpha_n}}|z|^{2\nu_n}G_{0,n}^{n,0}
\Big(\begin{aligned}& \rule{1cm}{0.1pt}\\-\nu_n,\nu_2-\nu_n&,....,\nu_{n-1}-\nu_n,0\end{aligned}\Big|\,|z|^2\Big),$$

\end{proposition}
 {\bf Proof.}
 Let $$\psi_{n}(x)=\int_0^{+\infty}e^{-\sum_{j=1}^{n-1}r_j^2-\frac{x^2}{r_1^2\cdots r_{n-1}^2}}\prod_{j=1}^{n-1}r_j^{2\nu_j-2\nu_n-1}dr_1\cdots dr_{n-1}.$$
Hence the Mellin transform
$$M_n(s)=\int_0^{+\infty}\psi_n(x)x^{s-1}ds,$$
is given by
$$M_{n}(s)=\int_0^{+\infty}\int_0^{+\infty}e^{-\sum_{j=1}^{n-1}r_j^2-\frac{x^2}{r_1^2\cdots r_{n-1}^2}}\prod_{j=1}^{n-1}r_j^{2\nu_j-2\nu_n-1}x^{s-1}dr_1\cdots dr_{n-1}dx.$$
Thus
$$M_{n}(s)=\int_0^{+\infty}\int_0^{+\infty}e^{-\sum_{j=1}^{n-1}r_j^2-u^2}\prod_{j=1}^{n-1}r_j^{2\nu_j-2\nu_n+s-1}u^{s-1}dr_1\cdots dr_{n-1}dx.$$
$$M_{n}(s)=\frac1{2^n}\Gamma(\frac s2)\prod_{j=1}^{n-1}\Gamma(\nu_j-\nu_n+\frac s2)=\frac1{2^n}\prod_{j=1}^{n}\Gamma(\nu_j-\nu_n+\frac s2).$$
Since, by the inverse Mellin transform one gets
$$\psi_n(x)=\frac1{2^{n-1}}G_{0,n}^{n,0}
\Big(\begin{aligned}& \rule{1cm}{0.1pt}\\-\nu_n,\nu_2-\nu_n&,....,\nu_{n-1}-\nu_n,0\end{aligned}\Big|\,x^2\Big).$$
The previous lemma gives the desired result.\\

We compute here the leading order asymptotic behaviour
of the weights $\omega_{n,N}(z)$ for $N$ large enough and for $z$ of order $O(1)$ using the saddle point method, such result will be used in the next section in order to study the limit of the one point correlation function. This method of proof can also be found in
Theorem 2 in reference [27] but we give it here for completeness.
\begin{lemma}\label{l1}
For large enough $N$, for every $z\in\Bbb C$ and fixed $\nu_2,...,\nu_n$, one has the following asymptotic relation
 $$\omega_n(z)\sim\frac{(2\pi)^{\frac 32(n-1)}}{2^{n-1}N^{\frac{n-1}2}\sqrt n}\frac{|z|^{\frac{1-n+2\alpha_n}{n}}}{\Big(1+|z|^{\frac2n}\Big)^{nN+\alpha_n+1}}.$$
where $\alpha_n=\sum\limits_{i=1}^n\nu_i$.
\end{lemma}
{\bf Proof.}\\
 We use the multidimensional representation of the weight, equation (\ref{eq4}).
$$\omega_{n,N}(z)=(2\pi)^{n-1}\int_0^{+\infty}dr_1\cdots dr_{n-1}\frac{|z|^{2\nu_n}}{(1+\frac{|z|^2}{r_1^2\cdots r_{n-1}^2})^{\nu_n+N+1}}\prod_{j=1}^{n-1}\frac{r_j^{2\nu_j-2\nu_n-1}}{(1+r^2_j)^{\nu_j+N+1}}.$$
It read to $$\omega_{n,N}(z)=(2\pi)^{n-1}\int_0^{+\infty}dr_1\cdots dr_{n-1}\frac{|z|^{2\nu_n}}{(1+\frac{|z|^2}{r_1^2\cdots r_{n-1}^2})^{\nu_n+1}}\prod_{j=1}^{n-1}\frac{r_j^{2\nu_j-2\nu_n-1}}{(1+r^2_j)^{\nu_j+1}}e^{-NS(r_1,...,r_{n-1})},$$
with $\displaystyle S=\log(1+\frac{|z|^2}{r_1^2\cdots r_{n-1}^2})+\sum_{j=1}^{n-1}\log(1+r_j^2).$

Now, we use the saddle point method.
A simple computation gives, for $j=1,...,n-1$
$$\frac{\partial S}{\partial r_j}=\frac{2r_j}{1+r_j^2}-\frac{2|z|^2}{r_j(r_1^2\cdots r_{n-1}^2+|z|^2)}$$
By symmetry one find the following solution to the equation $ \frac{\partial S}{\partial r_j}=0$, $r_j=r_*=|z|^{\frac 1n}$.
and the result follows.
For the details calculation of the saddle point, see appendix A below.

\section{Orthogonal polynomials}
From the fact that our weight $\omega_{n,N}(z)$ is angle-independent it immediately follows that the corresponding
orthogonal polynomials are monic $p_k(z) = z^k$. Indeed, the integral
$$\int_\C p_k(z)p_\ell(z^*)\omega_n(z)d^2z=\int_0^{2\pi}e^{i(k-\ell)\theta}d\theta\int_0^{+\infty}r\omega_n(r)dr=h_{k,n}\delta_{k \ell},$$
gives for $k\neq\ell$. One gets as in the previous the square norm,
$$h_{k,n}=\pi^n\big(\Gamma(N-k)\big)^n\prod_{i=1}^n\frac{\Gamma(\nu_i+k+1)}{\Gamma(\nu_i+N+1)}.$$
So that, the kernel $K_{n, N}(z_i,z_j)$ of orthonormal polynomials reads
\begin{equation}\label{eq6}K_{n,N}(z_i,z_j)=\frac1{\pi^n}\sqrt{\omega_n(z_i)\omega_n(z_j)}
\sum_{k=0}^{N-1}\frac1{\big(\Gamma(N-k)\big)^n}\prod_{i=1}^n\frac{\Gamma(\nu_i+N+1)}{\Gamma(\nu_i+k+1)}(z_iz_j^*)^k.\end{equation}
Recall that the generalized hypergeometric function is defined by
$${}_pF_q(a_1,...,a_p;b_1,...,b_q;z)=\sum_{k=0}^{+\infty}\frac{(a_1)_k...(a_p)_k}{(b_1)_k...(b_q)_k}\frac{z^k}{k!},$$
where $(a)_k=\frac{\Gamma(a+k)}{\Gamma(a)}=a(a+1)...(a+k-1)$ is the Pochammer symbol. If $a_1=...=a_p$ we use for convenient the notation for the hypergeometric function
${}_pF_q(a_1,...,a_p;b_1,...,b_q;x)={}_pF_q(a_1;b_1,...,b_q;x)$.
\begin{proposition}
For every $N\in\N$, $z_i,z_j\in\C$,
\small{$$K_{n,N}(z_i,z_j)=\frac{\prod_{i=1}^n(\nu_i+1)_N}{\pi^n\Gamma(N)^n}\sqrt{\omega_n(z_i)\omega_n(z_j)}\;\; {}_nF_{n-1}\big(1-N;1+\nu_2,...,1+\nu_n;(-1)^nz_iz_j^*\big).$$}
\end{proposition}
{\bf Proof.} First remark that for $a\in\C$, $(-a)_k=(-1)^ka(a-1)...(a-k+1)$. Since
$$\frac{\Gamma(N)}{\Gamma(N-k)}=(N-1)(N-2)...(N-k)=(-1)^k(1-N)_k,$$
$$\Gamma(1+\nu_1+k)=k!$$
and $$\frac{\Gamma(1+\nu_i)}{\Gamma(1+\nu_i+k)}=\frac1{(1+\nu_i)_k},$$
which give the desired result.

Remark that for $p=1,q=0$, $\displaystyle {}_1F_0(a;\,;z)=\frac1{(1-z)^a}$, moreover
$\displaystyle \omega_{1,N}(z)=\frac1{(1+|z|^2)^{N+1}}$, hence
$$K_{1, N}(z,z)=\frac{N}{\pi}{}_1F_0(1-N;\,;-|z|^2)=\frac{N}\pi\frac1{(1+|z|^2)^2}.$$
The one point correlation function or the eigenvalues density at finite $N$ read to
$$\rho_{1, 1, N}(z)=\frac1NK_{1, N}(z,z),$$ hence
$$\rho_{1, 1, N}(z)=\frac1{\pi(1+|z|^2)^2}.$$
The density of eigenvalues is independent of $N$.

In the general case the empirical density of eigenvalues is
\small{$$\rho_{n, N}(z):=\frac1N\rho_{n,1, N}(z)=\frac{\prod_{i=1}^n(\nu_i+1)_N}{\pi^nN\Gamma(N)^n}
\omega_{n,N}(z){}_nF_{n-1}\big(1-N;1+\nu_2,...,1+\nu_n;(-1)^n|z|^2\big).$$}
 normilized with the condition $\displaystyle\int_{\C}\rho_{n,N}(z)d^2z=1$.
\begin{proposition}\label{p1} Let n be a positive integer and $\nu_i=\lambda_iN$ with $\alpha_i\geq 0$, then for every $z\in\C$,
$$\lim_{N\to +\infty}\rho_{n,N}(z)=\rho_{n,\lambda}(z),$$
where the macroscopic density $\rho_{n,\alpha}$ is defined through its Mellin transform, for all $s\in]-\infty,\frac1n[$,
$$\int_\C|z|^{2s}\rho_{n,\lambda}(z)d^2z=\int_0^1\prod_{i=1}^n\frac{(t+\lambda_i)^s}{(1-t)^s}dt.$$
In particular, for $\lambda_1 = \cdots = \lambda_n = 0$ (this includes a product of square matrices), we
have
$$\rho_{n,0}(z)=\frac1{\pi n}\frac{|z|^{\frac{2-2n}{n}}}{\Big(1+|z|^{\frac 2n}\Big)^2}.$$
\end{proposition}
For $n=1$, one recovers the density of eigenvalues of the Cauchy-Lorentz ensemble for one random matrix.

{\bf Proof.} We saw from equation (\ref{eq6}) that
$$\rho_{n, N}(z):=\frac1NK_{n,N}(z,z)=\frac1{N\pi^n}\omega_{n,N}(z)
\sum_{k=0}^{N-1}\frac1{\big(\Gamma(N-k)\big)^n}\prod_{i=1}^n\frac{\Gamma(\nu_i+N+1)}{\Gamma(\nu_i+k+1)}|z|^{2k}.$$
The density is invariant under rotations in the complex, thus one can compute the Mellin transform of the radial density
$\varphi_{n,N}(r)=2\pi\displaystyle r\rho_{n, N}(r)$

 As in the previous one proves that for $s\in\C$ with $Re(s)<\frac1n$
$$\int_0^{+\infty}r^{2s}\varphi_{n, N}(r)dr=\frac1{N}
\sum_{k=0}^{N-1}\big(\frac{\Gamma(N-k-s)}{\Gamma(N-k)}\big)^n\prod_{i=1}^n\frac{\Gamma(\nu_i+k+s+1)}{\Gamma(\nu_i+k+1)}.$$
We would like to approximate the sum with an integral. Sine
the function $$t\mapsto \frac{\Gamma(N-t-s)}{\Gamma(N-t)}\prod_{i=1}^n\frac{\Gamma(\nu_i+t+s+1)}{\Gamma(\nu_i+t+1)},$$
is a non-negative, and monotonously increasing for $t\geq 0$, $N\in\N$ and $s\in]-\infty,\frac1n[$. Thus by
\begingroup\makeatletter\def\f@size{10}\check@mathfonts$$ \frac1N\int_0^N\frac{\Gamma(N-t-s+1)}{\Gamma(N-t+1)}\prod_{i=1}^n\frac{\Gamma(\nu_i+t+s)}{\Gamma(\nu_i+t)}dt\leq\int_0^{+\infty}r^{2s}\varphi_{n, N}(r)dr\leq \frac1N\int_0^N\frac{\Gamma(N-t-s)}{\Gamma(N-t)}\prod_{i=1}^n\frac{\Gamma(\nu_i+t+s+1)}{\Gamma(\nu_i+t+1)}dt.$$\endgroup
and by making the change of variable $t\rightarrow Nt$, and putting $\nu_i=\alpha_i N$, one gets
\begingroup\makeatletter\def\f@size{10}\check@mathfonts$$\int_0^1\frac{\Gamma(N(1-t)-s+1)}{\Gamma(N(1-t)+1)}
\prod_{i=1}^n\frac{\Gamma(N(\lambda_i+t)+s)}{\Gamma(N(\lambda_i+t))}dt\leq\int_{\C}r^{2s}\varphi_{n, N}(r)dr\leq \int_0^1\frac{\Gamma(N(1-t)-s)}{\Gamma(N(1-t))}\prod_{i=1}^n\frac{\Gamma(N(\lambda_i+t)+s+1)}{\Gamma(N(\lambda_i+t)+1)}dt.$$\endgroup
Using the following asymptotic formula for the gamma function $\displaystyle\Gamma(a+N)\sim N^a\Gamma(N),$ one gets for all $s\in]-\infty,\frac1n[$,
$$\lim_{N\to +\infty}\int_0^{+\infty}r^{2s}\varphi_{n,N}(r)dr=\int_0^1\prod_{i=1}^n\frac{(t+\lambda_i)^s}{(1-t)^s}dt.$$
Moreover the two functions in the left and right said are extend as a holomorphic functions on the half plane $Re(s)<\frac1n$, hence they are equals.
The macroscopic density is obtained by an inverse Mellin transform.

For $\lambda_1=...=\lambda_n=0$, one gets
$$\lim_{N\to +\infty}\int_0^{+\infty}r^{2s}\varphi_{n,N}(r)dr=\int_0^1\frac{t^{ns}}{(1-t)^{ns}}dt,$$
or $$\lim_{N\to +\infty}\int_\C|z|^{2s}\rho_{n,N}(z)d^2z=\int_0^1\frac{t^{ns}}{(1-t)^{ns}}dt,$$
Furthermore by making a change of variable to polar coordinates
$$\frac1{\pi n}\int_{\C}|z|^{2s}\frac{|z|^{\frac{2-2n}{n}}}{\Big(1+|z|^{\frac 2n}\Big)^2}d^2z=\frac2n\int_0^{+\infty}\frac{r^{2s+\frac2n-1}}{\Big(1+r^{\frac 2n}\Big)^2}dr,$$
subsituting $r=\frac{u^{\frac n2}}{(1-u)^{\frac n2}},$ one obtains
$$\frac1{\pi n}\int_{\C}|z|^{2s}\frac{|z|^{\frac{2-2n}{n}}}{\Big(1+|z|^{\frac 2n}\Big)^2}d^2z=\int_0^1\frac{u^{ns}}{(1-u)^{ns}}du.$$
Which prove that $\rho_{n,N}$ converges to the density $\displaystyle\frac1{\pi n}\frac{|z|^{\frac{2-2n}{n}}}{\Big(1+|z|^{\frac 2n}\Big)^2}$ in all the complex plane.

\section{Asymptotic of correlation functions}
Let $z_i=\big(\frac{\xi_i}{\sqrt {nN}}\big)^n$, then $d^2z_i=\frac{n^2}{(nN)^n}|\xi_i|^{2n-2}d^2\xi_i,$ and
$$\widehat \rho_{n,k,N}(\xi_1,...,\xi_k)=\frac{n^{2k}}{(nN)^{nk}}\prod_{i=1}^k|\xi_i|^{2n-2}\rho_{n,k,N}(\big(\frac{\xi_1}{\sqrt {nN}}\big)^n,...,\big(\frac{\xi_k}{\sqrt {nN}}\big)^n).$$
 we define $$\widehat K_{n,N}(\xi_i,\xi_j)=\frac{n^{2}}{(nN)^{n}}|\xi_i\xi_j|^{n-1}K_{n,N}(\big(\frac{\xi_i}{\sqrt {nN}}\big)^n,\big(\frac{\xi_i}{\sqrt {nN}}\big)^n).$$
 Then $$\widehat \rho_{n,k,N}(\xi_1,...,\xi_k)=\det\big(\widehat K_{n,N}(\xi_i,\xi_j)\big)_{i,j=1}^k.$$
\subsection{bulk limit}
We begin by the following lemma which will be useful in the next proposition. Let
$$\displaystyle f_{n,N}(x)=\sum_{k=0}^{N-1}\prod_{i=1}^n\frac{\Gamma(\nu_i+N+1)}{\Gamma(N-k)\Gamma(\nu_i+k+1)}x^{k}.$$
\begin{lemma} For fixed $\nu_i$, $i=2,...,n$ and for $N$ and $x\geq 0$ large enough we have
$$f_{n,N}(x)\sim \Big(\prod_{i=1}^n\frac{\Gamma(\nu_i+N+1)}{2\pi}\Big)\frac{\sqrt{2\pi}e^{n(N-1)}}{\sqrt n(N-1)^{nN+\alpha_n-\frac12}}x^{\frac{1-n-2\alpha_n}{2n}}(1+x^{\frac1n})^{nN+\alpha_n-1}.$$
\end{lemma}
{\bf Proof.}
We saw that
$$\Gamma(\nu_i+k+1)=(\nu_i+k)...(k+1)\Gamma(k+1),$$
Using Stirling formula, then for fixed $\nu_i$, and large enough $N$ and $k$  $$\Gamma(\nu_i+k+1)\sim k^{\nu_i}\sqrt{2\pi k}(k/e)^k,$$
and $$\Gamma(N-k)\sim\sqrt{2\pi (N-k-1)}\big((N-k-1)/e\big)^{N-k-1},$$
Since the serie $f_{n,N}(x)$ diverges (see appendix B), one can
approximating the sum by an integral and apply the saddle point method, hence:
$$\begin{aligned}f_{n,N}(x)&\sim \prod_{i=1}^n\frac{\Gamma(\nu_i+N)}{2\pi}\int_1^N\frac{x^k}{k^{\alpha_n}(k(N-k-1))^{\frac n2}(\frac ke)^{nk}(\frac{N-k-1}e)^{n(N-k-1)}}dk\\&=\prod_{i=1}^n\frac{\Gamma(\nu_i+N)}{2\pi}\int_1^N\frac{e^{S(k)}}{k^{\alpha_n}(k(N-k-1))^{\frac n2}}dk,\end{aligned}$$
with $\displaystyle\alpha_n=\sum_{i=1}^n\nu_i$, and $$\displaystyle S(k)=-nk\log k-n(N-k-1)\log(N-k-1)+nk+n(N-k-1)+k\log x.$$
Since $$S'(k)=-n\log k+n\log(N-k-1)+\log x,$$
hence $S'(k)=0$ iff $\displaystyle k^*:=k=(N-1)\frac{x^{\frac1n}}{1+x^{\frac1n}}.$ Moreover
$$e^{S(k^*)}=\frac{e^{n(N-1)}}{(N-1)^{n(N-1)}}(1+x^{\frac1n})^{n(N-1)}.$$
and
$$S''(k^*)=-\frac n{N-1}\frac{(1+x^{\frac 1n})^2}{x^{\frac 1n}}.$$
By a simple computation one has
$$\frac{1}{k^{\alpha_n}(k(N-k-1))^{\frac n2}}=\frac1{(N-1)^{n+\alpha_n}}\frac{(1+x^{\frac1n})^{\alpha_n+n}}{ x^{\frac{\alpha_n}n+\frac12}}.$$
Putting all this together and using the saddle point asymptotic expression one gets
$$f_{n,N}(x)\sim \Big(\prod_{i=1}^n\frac{\Gamma(\nu_i+N+1)}{2\pi}\Big)\frac{\sqrt{2\pi}e^{n(N-1)}}{\sqrt n(N-1)^{nN+\alpha_n-\frac12}}x^{\frac{1-n-2\alpha_n}{2n}}(1+x^{\frac1n})^{nN+\alpha_n-1}.$$
This complete the proof.
\begin{proposition}\label{p2}
$$\lim_{N\to+\infty}N\widehat K_{n,N}(\xi_i,\xi_j)=\frac1\pi\frac{(\xi_i\xi_j^*)^{\frac{1-n-2\alpha_n}2}}{|\xi_i\xi_j|^{\frac{1-n-2\alpha_n}2}}e^{-\frac12\big(|\xi_i|^2+|\xi_j|^2-2\xi_i\xi_j^*\big)}.$$
And $$\lim_{N\to +\infty}N^k\widehat \rho_{n,k,N}(\xi_1,...,\xi_k)=\det\Big(\K(\xi_i,\xi_j)\Big)_{ i,j=1}^k,$$
where $\K(x,y)=\frac1\pi e^{-\frac12\big(|x|^2+|y|^2-2xy^*\big)}.$
\end{proposition}
\begin{remark} One can see the two point correlation function is
$$\rho_2(\xi_1,\xi_2):=\lim_{N\to +\infty}N^2\widehat \rho_{n,k,N}(\xi_1,\xi_2)=\frac1{\pi^2}\Big(1-e^{-|\xi_1-\xi_2|^2}\Big),$$
which describes correlations of eigenvalues at distances $|\xi_1-\xi_2|$ of order unity.
\end{remark}
\begin{remark}
The kernel is unitary equivalent to the Ginibre kernel. When one calculates the correlation functions
$\widehat\rho_{n,k}(\xi_1, . . . , \xi_k)$ all phase factors cancel and one obtains exactly the same k-point correlation
functions as for the Ginibre ensemble with $n =1$. Therefore all correlation functions of the product
matrix $Y_n$ in the bulk limit are universal.
\end{remark}
{\bf Proof of proposition \ref{p2}}. We have seen, the $k$-correlation function is determinental, that is
$$\rho_{n,k,N}(z_1,...,z_k)=\det\big(K_{n,N}(z_i,z_j)_{i,j=1}^k\big),$$
where $$K_{n,N}(z_i,z_j)=\frac1{\pi^n}\sqrt{\omega_n(z_i)\omega_n(z_j)}\sum_{k=0}^{N-1}\frac1{\big(\Gamma(N-k)\big)^n}
\prod_{i=1}^n\frac{\Gamma(\nu_i+N+1)}{\Gamma(\nu_i+k+1)}
(z_iz_j^*)^k.$$
Hence $$K_{n,N}(z_i,z_j)=\frac{1}{\pi^n}\sqrt{\omega_n(z_i)\omega_n(z_j)}f_{n,N}(z_iz_j^*).$$
We saw from lemma \ref{l1} that when $N\to+\infty$, and for every fixed $z\in\C$, $$\omega_{n,N}(z)\sim a_{n,N}\frac{|z|^{\frac{1-n+2\alpha_n}n}}{\Big(1+|z|^{\frac2n}\Big)^{nN+\alpha_n+1}},$$
and from the previous proposition one has $$f_{n,N}(z)\sim b_{n,N}|z|^{\frac{1-n-2\alpha_n}{n}}(1+|z|^{\frac2n})^{nN+\alpha_n-1}$$
with $\displaystyle a_{n,N}=\frac{(2\pi)^{\frac 32(n-1)}}{2^{n-1}N^{\frac{n-1}2}\sqrt n}$ and $\displaystyle b_{n,N}=\frac{\prod_{i=1}^n\Gamma(\nu_i+N+1)}{(2\pi)^{n-\frac 12}\sqrt n}\frac{e^{n(N-1)}}{(N-1)^{nN-\frac12}}$. Thus, it is enough to study the asymptotic of the function
$$g_{n,N}(z_i,z_j)=\sqrt{\frac{|z_i|^{\frac{1-n+2\alpha_n}n}}{\Big(1+|z_i|^{\frac2n}\Big)^{n(N-1)+\alpha_n+1}}
\frac{|z_j|^{\frac{1-n+2\alpha_n}n}}{\Big(1+|z_j|^{\frac2n}\Big)^{n(N-1)+\alpha_n+1}}}(z_iz_j^*)^{\frac{1-n-2\alpha_n}{2n}}
(1+(z_iz_j^*)^{\frac1n})^{n(N-1)+\alpha_n-1}.$$
Let $\displaystyle z=\big(\frac{\xi}{\sqrt{nN}}\big)^n,$
then we obtain
$$g_{n,N}(\xi_i,\xi_j)=\frac1{(nN)^{1-n}}\sqrt{\frac{|\xi_i\xi_j|^{1-n+2\alpha_n}}
{\Big(1+\frac{|\xi_i|^2}{nN}\Big)^{n(N-1)+\alpha_n+1}\Big(1+\frac{|\xi_j|^2}{nN}\Big)^{n(N-1)+\alpha_n+1}}}
(\xi_i\xi_j^*)^{\frac{1-n-2\alpha_n}{2}}\Big(1+\frac{\xi_i\xi_j^*}{nN}\Big)^{n(N-1)+\alpha_n-1}.$$
As $N\to +\infty$, since $\alpha_n$ is fixed, it follows that
$$g_{n,N}(\xi_i,\xi_j)\sim\frac1{(nN)^{1-n}}|\xi_i\xi_j|^{\frac{1-n+2\alpha_n}2}
(\xi_i\xi_j^*)^{\frac{1-n-2\alpha_n}2}e^{-\frac12\big(|\xi_i|^2+|\xi_j|^2-2\xi_i\xi_j^*\big)}.$$
Moreover $$N\widehat K_{n,N}(\xi_i,\xi_j)=Na_{n,N}b_{n,N}\frac{n^{2}}{(nN)^{n}}|\xi_i\xi_j|^{n-1}g_{n,N}(\xi_i,\xi_j).$$
Using again the Stiriling formula, it yields $$\prod_{i=1}^n\frac{\Gamma(\nu_i+N+1)}{\sqrt{2\pi}}=\Big(\frac{\Gamma(N)}{\sqrt{2\pi}}\Big)^n\prod_{i=1}^n(\nu_i+N)(\nu_i+N-1)...N\sim N^{\alpha_n+n}\frac{(N-1)^{\frac n2}(N-1)^{n(N-1)}}{e^{n(N-1)}},$$
Substituting in $a_{n,N}$ and $b_{n,N}$, one has the following asymptotic formula $$\displaystyle\frac{1}{\pi^n} a_{n,N}b_{n,N}\sim\frac1{n\pi}.$$ Which gives as $N$ goes to infinity
$$N\widehat K_{n,N}(\xi_i,\xi_j)\sim\frac1\pi\frac{(\xi_i\xi_j^*)^{\frac{1-n-2\alpha_n}2}}{|\xi_i\xi_j|^{\frac{1-n-2\alpha_n}2}}e^{-\frac12\big(|\xi_i|^2+|\xi_j|^2-2\xi_i\xi_j^*\big)}.$$
Since the angular phase disappear in the determinant formula. In fact it is a multiplication by diagonal unitary transformation. This gives the desired result.
\subsection{Limit at the origin}
\begin{proposition}
For $\xi_i,\xi_{j}$ in compact set of $\Bbb C$, and every  fixed $\nu_2,...,\nu_n$
$$\K_n(\xi_i,\xi_j)=\lim_{N\to+\infty}K_{n,N}\big(\frac{\xi_i}{N^{\frac n2}},\frac{\xi_j}{N^{\frac n2}}\big)=\frac{\sqrt{\omega_n(\xi_i)\omega_n(\xi_j) }}{\pi\prod_{i=1}^n\Gamma(1+\nu_i)}\;{}_0F_{n-1}\Big(-;1+\nu_2,...,1+\nu_n;\xi_i\xi_j^*\Big),$$
where $\omega_n(\xi)=|\xi|^{2\nu_n}G_{0,n}^{n,0}
\Big(\begin{aligned}& \rule{1cm}{0.1pt}\\-\nu_n,\nu_2-\nu_n&,....,\nu_{n-1}-\nu_n,0\end{aligned}\Big|\,|\xi|^2\Big).$
\end{proposition}
{\bf Proof.} Let $z=\frac \xi{N^{\frac n2}}$, then $d^2z=\frac1{N^ n}d^2\xi$.

We saw in proposition (\ref{pr01}) that, as $N$ goes to $+\infty$
\begin{equation}\label{eq1}\omega_{n,N}\big(\frac z{N^{\frac n2}}\big)\sim \frac{\pi^{n-1}}{N^{\alpha_n}}|z|^{2\nu_n}G_{0,n}^{n,0}
\Big(\begin{aligned}& \rule{1cm}{0.1pt}\\-\nu_n,\nu_2-\nu_n&,....,\nu_{n-1}-\nu_n,0\end{aligned}\Big|\,|z|^2\Big).\end{equation}
Furthermore
 \small{\begin{equation}\label{eq2}K_{n,N}\big(\frac{\xi_i}{N^{\frac n2}},\frac{\xi_j}{N^{\frac n2}}\big)=d_{n,N}\sqrt{\omega_{n,N}(\frac{\xi_i}{N^{\frac n2}})\omega_{n,N}(\frac{\xi_j}{N^{\frac n2}})}\;\; {}_nF_{n-1}\big(1-N;1+\nu_2,...,1+\nu_n;(-1)^n\frac{\xi_i\xi_j^*}{N^{ n}}\big),\end{equation}}
 with $\displaystyle d_{n,N}=\frac{\prod_{i=1}^n\Gamma(\nu_i+N+1)}{\pi^nN^n\prod_{i=2}^n\Gamma(1+\nu_i)\Gamma(N)^n}$.

Recall that $${}_nF_{n-1}\big(1-N;1+\nu_2,...,1+\nu_n;(-1)^n\frac{x}{N^{ n}}\big)=\sum_{k=0}^{N-1}\frac{\big(( 1-\frac2N)...(1-\frac{k}{N})\big)^n}{\prod_{i=2}^n(1+\nu_i)_k}\frac{x^k}{k!}.$$
 By the equation below, valid for all $N$ and $k\leq N-1$, $$( 1-\frac2N)...(1-\frac{k}{N})\leq 1,$$
 we have on every compact of $\C$ \begin{equation}\label{eq3}\lim_{N\to+\infty}{}_nF_{n-1}\big(1-N;1+\nu_2,...,1+\nu_n;(-1)^n\frac{x}{N^{ n}}\big)={}_0F_{n-1}\big(-\;;1+\nu_2,...,1+\nu_n;x\big).\end{equation}
 Furthermore \begin{equation}\label{eq0}d_{n,N}=\frac{\prod_{i=1}^n\Gamma(\nu_i+N+1)}{\pi^nN^n\prod_{i=2}^n\Gamma(1+\nu_i)\Gamma(N)^n}\sim\frac{
  N^{\alpha_n+n}\Gamma(N)^n}{\pi^nN^n\prod_{i=2}^n\Gamma(1+\nu_i)\Gamma(N)^n}=\frac{N^{\alpha_N}}{\pi^n\prod_{i=2}^n\Gamma(1+\nu_i)},\end{equation}
From equations (\ref{eq1}), (\ref{eq2}), (\ref{eq2}) and (\ref{eq0}), one gets the desired result.
\begin{remark}
1) Remark that if all the matrices $X_i$ are square $\nu_i=0$, for all $i=1,..,n$, then,
the Kernel read to
$$\K_n(\xi_i,\xi_j)=\frac1{\pi}\sqrt{G_{0,n}^{n,0}
\Big(\begin{aligned}&\;-\;\\&\;0\;\end{aligned}\Big|\,|\xi_i|^2\Big)G_{0,n}^{n,0}
\Big(\begin{aligned}&\;-\;\\&\;0\;\end{aligned}\Big|\,|\xi_j|^2\Big)}\;{}_0F_{n-1}\Big(-;1;\xi_i\xi_j^*\Big),$$
2) For $n=1$ one gets the Ginibre Kernel $\K_1(\xi_i,\xi_j)=\frac1\pi e^{-\frac12\big(|\xi_i|^2+|\xi_j|^2-2\xi_i\xi^*_j\big)}$, in fact ${}_0F_0\Big(-;-;x\Big)=e^x$, and $G_{0,1}^{1,0}\Big(\begin{aligned}&-\\& 0\end{aligned}\;\;\Big|\,x\Big)=e^{-x}$.\\
3) For $n=2$, since ${}_0F_{1}\Big(-;1+\nu_2;x\Big)=\Gamma(1+\nu_2)(\frac{x}{2})^{-\nu_2}J_{\nu_2}(-2\sqrt x)$, and  $G_{0,2}^{2,0}
\Big(\begin{aligned}&\;\;\;\;-\;\\&\;-\nu_2\;0\end{aligned}\Big|\,x\Big)=2x^{-\frac{\nu_2}2}K_{-\nu_2}(2\sqrt x)$
and $$\K_2(\xi_i,\xi_j)=\frac{2}{\pi}\big(\frac{|\xi_i\xi_j|}{\xi_i\xi_j^*}\big)^{\frac{\nu_2}2}\sqrt{K_{-\nu_2}(2|\xi_i|)K_{-\nu_2}(2|\xi_j|)}
\;\;J_{\nu_2}(-2\sqrt{\xi_i\xi_j^*}).$$
moreover the correlation function read to,
$$\widehat\rho_{2,k}(\xi_1,...,\xi_k)=\det\Big(\widetilde\K_2(\xi_i,\xi_j)\Big)_{1\leq i,j\leq k},$$
with $\displaystyle \widetilde\K_2(\xi_i,\xi_j)=\sqrt{K_{-\nu_2}(2|\xi_i|)K_{-\nu_2}(2|\xi_j|)}
\;\;J_{\nu_2}(-2\sqrt{\xi_i\xi_j^*})$, in that case the correlation functions are expressed essentially in term of Bessel functions.
\end{remark}

Moreover $\K_1(\xi,\xi)=\frac1\pi$ which is the one point correlation function. If we look for the density of eigenvalues one has $\rho_1(0)=\frac1\pi$ see proposition \ref{p1}

\begin{center}Appendix A\end{center}
Evaluation of the asymptotic of $\omega_{n,N}(z)$ when $N$ goes to infinity.

We saw that,
$$\omega_{n,N}(z)=(2\pi)^{n-1}|z|^{2\nu_n}\int_0^{+\infty}dr_1\cdots dr_{n-1}\frac1{(1+\frac{|z|^2}{r_1^2...r_{n-1}^2})^{\nu_n+1}}\prod_{j=1}^{n-1}\frac{r_j^{2\nu_j-2\nu_n-1}}{(1+r^2_j)^{\nu_j+1}}e^{-NS(r_1,...,r_{n-1})},$$
where $\displaystyle S(r_1,...,r_{n-1})=\sum_{j=1}^{n-1}\log(1+r_j^2)+\log (1+\frac{|z|^2}{r_1^2...r_{n-1}^2}).$
And the solution to the problem
$\displaystyle\frac{\partial S}{\partial r_j}=0$ is $r_j=r_*=|z|^{\frac 1n}$ for every $j$.
Hence, by the saddle point method one has
\begin{equation}\label{eq5}\omega_{n,N}(z)\sim(2\pi)^{n-1}\big(\frac{2\pi}N\big)^{\frac{n-1}{2}}\frac{|z|^{\frac{2\alpha_n-n+1}{n}}}{(1+|z|^{\frac2n})^{\alpha_n+n}}
\frac{e^{-NS(r_*,...,r_*)}}{\sqrt{\det(S'')(r_*,...,r_*)}}.\end{equation}
Moreover $S(r_*,...,r_*)=(n-1)\log(1+r_*^2)+\log (1+\frac{|z|^2}{r_*^{2n-2}}).$
This gives
$$e^{-NS(r_*,...,r_*)}=\frac1{(1+|z|^{\frac2n})^{nN}}.$$
We saw that
$$\frac{\partial S}{\partial r_j}=\frac{2r_j}{1+r_j^2}-\frac{2|z|^2}{r_j(r_1^2\cdots r_{n-1}^2+|z|^2)},$$
and for $i\neq j$,
$$\frac{\partial^2 S}{\partial r_i\partial r_j}=\frac{4|z|^2r_1^2\cdots r_{n-1}^2}{r_ir_j(r_1^2\cdots r_{n-1}^2+|z|^2)^2},$$
$$\frac{\partial^2 S}{\partial^2 r_j}=\frac{2-2r^2_j}{(1+r_j^2)^2}+\frac{4|z|^2r_1^2\cdots r_{n-1}^2}{r^2_j(r_1^2\cdots r_{n-1}^2+|z|^2)^2}+\frac{2|z|^2}{r^2_j(r_1^2\cdots r_{n-1}^2+|z|^2)},$$
hence
$$\frac{\partial^2 S}{\partial r_i\partial r_j}\mid_{r_j=|z|^{\frac1n}}=\frac4{(1+|z|^{\frac2n})^2},$$
$$\frac{\partial^2 S}{\partial^2 r_j}\mid_{r_j=|z|^{\frac1n}}=\frac8{(1+|z|^{\frac2n})^2},$$
$$\sqrt{\det(S'')}=\frac{2^{n-1}\sqrt n}{(1+|z|^{\frac2n})^{n-1}}.$$
One gets from equation (\ref{eq5}) as $N$ and $|z|$ large enough
$$\omega_{n,N}(z)\sim\frac{(2\pi)^{\frac32(n-1)}}{2^{n-1}{N^{\frac{n-1}2}\sqrt n}}
\frac{|z|^{\frac{2\alpha_n-n+1}{n}}}{(1+|z|^{\frac2n})^{nN+\alpha_n+1}}.$$
\begin{center}Appendix B\end{center}
Recall the series of $f_{n,N}$ is given by
$$\displaystyle f_{n,N}(x)=\sum_{k=0}^{N-1}\prod_{i=1}^n\frac{\Gamma(\nu_i+N+1)}{\Gamma(N-k)\Gamma(\nu_i+k+1)}x^{k}.$$
Moreover As $N$ and $k$ large enough
$$\frac1{\Gamma(\nu_i+k+1)}=\frac1{(\nu_i+k)...(k+1)\Gamma(k+1)}\sim \frac1{k^{\nu_i}k!},$$
$$\Gamma(\nu_i+N)\sim (N-1)^{\nu_i}(N-1)!,$$
and $$\frac1{\Gamma(N-k)}\sim \frac{N^{k}}{(N-1)!},$$
hence $$\prod_{i=1}^n\frac{\Gamma(\nu_i+N+1)}{\Gamma(N-k)\Gamma(\nu_i+k+1)}\sim\frac{N^{\alpha_n+n}N^{nk}}{k^{\alpha_n}(k!)^n}.$$
Furthermore for every $k\leq N-1$ and all $n\geq 1$,
$$\frac{(N-1)^{\alpha_n}N^{n(k+1)}}{k^{\alpha_n}(k!)^n}\geq N^n\frac{ 1}{k!},$$
It follows that for every $x\geq 0$, and for large $N$
$$f_{n,N}(x)\geq N^n\sum_{k=0}^{N-1}\frac{ x^k}{k!}.$$
And $$\lim_{N\to +\infty}f_{n,N}(x)=+\infty,\qquad\forall\;x\geq 0.$$

Address:  College of Applied Sciences
   Umm Al-Qura University
  P.O Box  (715), Makkah,
  Saudi Arabia.\\
  Facult\'e des Scinces de Tunis,  Campus Universitaire El-Manar, 2092 El Manar Tunis.\\
E-mail: bouali25@laposte.net \& mabouali@uqu.edu.sa

\end{document}